\documentclass[12pt,leqno,a4paper]{article}

\usepackage{latexsym, amsbsy}

\usepackage{amssymb}
\usepackage{amsmath}
\usepackage{graphicx}
\usepackage{amsfonts}
\newtheorem{theo}{Theorem}
\newtheorem{obs}{Observation}
\newtheorem{defin}{Definition}
\newtheorem{lema}{Lemma}
\newtheorem{con}{Consequence}
\newtheorem{cor}{Corollary}
\newtheorem{prop}{Proposition}
\newtheorem{rema}{Remark}
\newtheorem{exem}{Example}

\newcommand{\ble}{\begin{lema}$\!\!\!\textrm{\bf{.}}$}
\newcommand{\ele}{\end{lema}}
\newcommand{\bde}{\begin{defin}$\!\!\!\textrm{\bf{.}}$}
\newcommand{\ede}{\end{defin}}
\newcommand{\bte}{\begin{theo}$\!\!\!\textrm{\bf{.}}$}
\newcommand{\ete}{\end{theo}}
\newcommand{\bob}{\begin{obs}$\!\!\!\textrm{\bf{.}}$}
\newcommand{\eob}{\end{obs}}
\newcommand{\bco}{\begin{cor}$\!\!\!\textrm{\bf{.}}$}
\newcommand{\eco}{\end{cor}}
\newcommand{\bcon}{\begin{con}$\!\!\!\textrm{\bf{.}}$}
\newcommand{\econ}{\end{con}}
\newcommand{\bre}{\begin{rema}$\!\!\!\textrm{\bf{.}}$}
\newcommand{\ere}{\end{rema}}
\newcommand{\bpr}{\begin{prop}$\!\!\!\textrm{\bf{.}}$}
\newcommand{\epr}{\end{prop}}
\newcommand{\bex}{\begin{exem}$\!\!\!\textrm{\bf{.}}$}
\newcommand{\eex}{\end{exem}}

\newcommand{\ds}{\displaystyle}

\date{}
\begin{document}

\title{On the Constant in The Lower Estimate for the Bernstein Operator}
\author{ Sorin G. Gal and Gancho T. Tachev}
\date{}
\maketitle

\begin{abstract}
For functions belonging to the classes $C^{2}[0, 1]$ and $C^{3}[0, 1]$, we establish the lower estimate with an explicit constant in approximation by Bernstein polynomials in terms of the second order Ditzian-Totik modulus of smoothness. Several applications to some concrete examples of functions are presented.
\end{abstract}

41A25, 41A36.


\section{Introduction}\medskip

For every function $f\in C[0,1]$ the Bernstein polynomial operator
is given by
$$ B_n(f;x)=\sum\limits_{k=0}^nf(\frac{k}{n})\cdot {{n}\choose{k}}x^k(1-x)^{n-k},\,x\in[0,1]\eqno(1.1)$$
To handle this operator, it is useful to utilize the second order Ditzian-Totik modulus of smoothness $\omega^2_{\varphi}(f,t)$ which is equivalent to the $K-$ functional given by
$$
K^2_{\varphi}(f,t^2):=\inf\limits_{g\in W^2_{\infty}(\varphi)}\left\{\|f-g\|+t^2\|\varphi^2g''\|_{L_{\infty}[0,1]}\right\},\eqno(1.2)$$
where $\|\cdot \|$ denotes the uniform norm on $C[0, 1]$, $\varphi(x)=\sqrt{x(1-x)},\,t>0, n\in N$ and the weighted Sobolev space $W^2_{\infty}(\varphi)$, is given by
$$W^2_{\infty}(\varphi)=\left\{f\in AC[0,1],\, f'\in AC_{loc},\, \varphi^2f''\in L_{\infty}[0,1]\right\}.\eqno(1.3)$$

{\bf Remark 1.} For definition, properties and many applications of $\omega^2_{\varphi}(f,t)$ and $K^2_{\varphi}(f,t^2)$ see [3]. For example, it is known that if in the definition of $K^2_{\varphi}(f,t^2)$, instead of the space $W^2_{\infty}(\varphi)$ we consider the space
$$C^{2}[0, 1]=\{g:[0, 1]\to \mathbb{R} ; g \mbox{ is twice continuously diferentiable on } [0, 1]\}$$
where $C^{2}[0, 1]\subset W^2_{\infty}(\varphi)$, this replacement does not have any effect on the equivalence between
$\omega^2_{\varphi}(f,t)$ and $K^2_{\varphi}(f,t^2)$.

In a pioneering work of Ditzian and Ivanov [2], a general theory was developed to obtain strong converse inequalities for a broad class of operators. For the Bernstein operator, in [2] it was proved a strong converse inequality of type B. The first proof of strong converse inequality of type A for Bernstein operator was given in 1994 by Knoop and Zhou in [15] and Totik in [19], which we cite here as :

{\bf Theorem A.} {\it There exist two absolute constants $C_{1}, C_{2}>0$ such that
$$C_{1}\omega^2_{\varphi}(f,\frac{1}{\sqrt{n}})\le \|f-B_nf\|\le C_{2}\omega^2_{\varphi}(f,\frac{1}{\sqrt{n}})\eqno(1.4)$$
holds for all $f\in C[0,1]$ and all} $n\in N$.\\[1ex]
The proof of Theorem A is very complicated. Concerning the absolute constants, very recently using different methods, in [17] among others it was proved that $C_{2}=3$ could be placed in the right-hand side of (1.4). But as far as we know, nothing is known about the constant in the lower estimate in (1.4). It is the aim of this paper to establish for the first time concrete value of the constant in the left-hand side of (1.4), but only for the functions $f\in C^{2}[0, 1]$ and $f\in C^{3}[0, 1]$.

Our main results can be stated as follows.

{\bf Theorem 1.} {\it For any $f\in C^{2}[0, 1]$ and any $\mu_{0}\in (0, 1)$, there exists $n_{1}(f, \mu_{0})\in \mathbb{N}$ (depending on $f$ and $\mu_{0}$), such that for all $n\ge n_1$, we have
$$\frac{\mu_{0}}{32}\cdot\omega^2_{\varphi}(f,\frac{1}{\sqrt{n}})\le \|f-B_nf\|\le 3\cdot\omega^2_{\varphi}(f,\frac{1}{\sqrt{n}}).\eqno(1.5)$$}

{\bf Corollary 1.} {\it For $m>0$, let us define the class of functions
$$C^{3, 2}_{M, m}[0, 1]=\{f:[0, 1]\to \mathbb{R} ; f\in C^{3}[0, 1], \|f^{\prime \prime \prime}\|\le M,\,|f^{\prime \prime}(x)|\ge m, x\in[0,1] \}.$$
Then for all $f\in C^{3, 2}_{M, m}[0, 1]$ and $n\ge n_{1}$ with $n_{1}=\left [\frac{1024 M^{2}}{m^{2}}\right ]+1$,
we have
$$\frac{1}{64}\cdot\omega^2_{\varphi}(f,\frac{1}{\sqrt{n}})\le \|f-B_nf\|\le 3\omega^2_{\varphi}(f,\frac{1}{\sqrt{n}}).$$}
{\bf Remark 2.} Evidently that Theorem 1 and Corollary 1 have the disadvantage that the constants $\frac{\mu_{0}}{32}$ and $\frac{1}{64}$ in the
corresponding lower estimates and the index $n_{1}$ depend on the functions $f$, as being valid only for $n\ge n_{1}$ and for functions in $C^{2}[0, 1]$ and in $C^{3}[0, 1]$, respectively.
However, their advantages is that it is for the first time when for a broad
class of functions (including many elementary concrete functions, like $exp(x)$, $sin(x)$, $cos(x)$,
$arctg(x)$ and so on..) these results allow to give concrete lower bounds in terms of the modulus $\omega^2_{\varphi}(f,\frac{1}{\sqrt{n}})$
for the norm of $\|B_nf-f\|$. Note that these constants are not possible to be deduced from all the other
well known results, like those from Knoop and Zhou \cite{knoop}, Totik \cite{to}, Ditzian and Ivanov \cite{zd1},
and many others.

In Section 2 we give some auxiliary results and establish in Theorems 3 and 4 norm estimates in Voronovskaja's theorem for Bernstein operator. In Section 3 we prove Theorem 1 and Corollary 1. Finally, in Section 4 several applications to some concrete examples of functions are presented.

\section{Auxiliary results}

As already mentioned, the moduli $\omega^2_{\varphi}(f,t)$ and $K^2_{\varphi}(f,t^2)$ are equivalent. For example Theorem 6.2 in Chapter 6 in [1] (see also [3]) states that, there are constants $C_1,C_2>0$, such that for all $f\in L_{\infty}$
$$C_1\omega^2_{\varphi}(f,t)\le K^2_{\varphi}(f,t^2) \le C_2\omega^2_{\varphi}(f,t),\,0<t\le \frac{1}{2}.\eqno(2.1)$$
For our goals it is important to determine explicitly the values of $C_1$ and $C_2$ in (2.1). We recall the following result, established in [17] and in [9] :

{\bf Theorem B.} {\it For all $f\in W^2_{\infty}(\varphi)$ and all $n\ge 1$ the following holds true:
$$A_1\omega^2_{\varphi}(f,\frac{1}{\sqrt{n}})\le K^2_{\varphi}(f,\frac{1}{n}) \le A_2\omega^2_{\varphi}(f,\frac{1}{\sqrt{n}}),\eqno(2.2)$$
where the value of $A_1=\frac{1}{16}$ was established in $[17]$ and $A_2=10$ follows from $[9]$.}

The next estimate was proved in Theorem 6.1 in [1] :

{\bf Theorem C.} {\it There is a constant $C>0$, depending only on $r$, such that for each $f\in W^2_p(\varphi),\,1\le p\le \infty$ we have }
$$\omega^r_{\varphi}(f,t)_{p}\le Ct^r\|\varphi^rf^{(r)}\|_p,\, 0\le t\le \frac{1}{2r}.\eqno(2.3)$$

Note that from the proof of Theorem C, it is not possible to determine the magnitude of $C$ in (2.3). But using Theorem B, we may proceed as follows : For $f\in C^{2}[0, 1]$ we have
$$\omega^2_{\varphi}(f,\frac{1}{\sqrt{n}})\le \frac{1}{A_1}\cdot K^2_{\varphi}(f,\frac{1}{n})\le \frac{1}{A_1}\cdot\frac{1}{n}\|\varphi^2f''\|.\eqno(2.4)$$
For the value of $A_1=\frac{1}{16}$ see Theorem 2 in [17] and for the value of $A_2=10$ see Corollary 9 in [9].

Even more, it is possible to establish a lower bound for $\omega^2_{\varphi}(f,\frac{1}{\sqrt{n}}),\, f\in C^{2}[0, 1]$, as follows.

{\bf Theorem 2.} {\it For any $f\in C^{2}[0, 1]$ and any $\lambda_{0}\in (0, 1)$, there exists $n_{0}$ (depending on $f$ and $\lambda_{0}$) such that for all $n\ge n_{0}$ the following
$$\lambda_{0}\cdot \frac{1}{n}\|\varphi^2f^{\prime \prime}\|\le \omega^2_{\varphi}(f,\frac{1}{\sqrt{n}})\eqno(2.5)$$
holds true.}

{\bf Proof:} If $f$ is a polynomial of degree $\le 1$ on $[0, 1]$, then the inequality one reduces to the equality $0=0$.

Therefore, suppose that $f$ is not a polynomial of degree $\le 1$. By the definition of $\omega^2_{\varphi}(f, \delta)$, for all $0\le \delta <1$ we can write (see e.g. \cite{D1})
$$\omega^2_{\varphi}(f, \delta)=\sup \{\sup  \{|f(x+h\varphi(x))-2f(x)+f(x-h\varphi(x))|; $$
$$x\in [h^{2}/(1+h^{2}), 1/(1+h^{2})] \} ; 0\le h \le \delta \}.$$
Because $f$ is not a polynomial of degree $\le 1$, it follows that $f^{\prime \prime}(x)$ is not identical equal to zero on $[0, 1]$ and that there exists a point $x_{0}\in [0, 1]$ such that $\|\varphi^{2} f^{\prime \prime}\|_{\infty}=|\varphi^{2}(x_{0})f^{\prime \prime}(x_{0})|>0$ (contrariwise would easily follow that $f$ is a polynomial of degree $\le 1$, a contradiction).

Since evidently that $x_{0}\in (0, 1)$, this implies that there exists $h_{0}\in (0, 1)$ such that for all
$h\in [0, h_{0}]$ we have $x_{0}\in [h^{2}/(1+h^{2}), 1/(1+h^{2})]$, and
$$\sup  \{|f(x+h\varphi(x))-2f(x)+f(x-h\varphi(x))|; x\in [h^{2}/(1+h^{2}), 1/(1+h^{2})] \}$$
$$\ge |f(x_{0}+h\varphi(x_{0})-2f(x_{0})+f(x-h\varphi(x_{0})|=h^{2}\varphi^{2}(x_{0})\cdot |f^{\prime \prime}(\xi_{h, x_{0}})|,$$
for all $h\in [0, h_{0}]$,
where from the mean value theorem
$$\xi_{h, x_{0}}\in [x_{0}-h\varphi(x_{0}), x+h\varphi(x_{0})].$$
For $h\to 0$, evidently that $\xi_{h, x_{0}}\to x_{0}$ and from the continuity of $f^{\prime \prime}$ on $[0, 1]$, it follows that
$\lim_{h\to 0}|f^{\prime \prime}(\xi_{h, x_{0}})|=|f^{\prime \prime}(x_{0}|>0$. Therefore, for $\lambda_{0}\in (0, 1)$,
there exists $0\le h_{1}< h_{0}$, such that for all $0\le h < h_{1}$, we have $|f^{\prime \prime}(\xi_{h, x_{0}})|\ge \lambda_{0}|f^{\prime \prime}(x_{0})|$ and combined with the above lower estimate, implies
$$\sup  \{|f(x+h\varphi(x))-2f(x)+f(x-h\varphi(x))|; x\in [h^{2}/(1+h^{2}), 1/(1+h^{2})] \}$$
$$\ge \lambda_{0}h^{2}\|\varphi^{2} f^{\prime \prime}\|.$$
Now, let $n_{0}\in \mathbb{N}$ be the smallest natural number such that $\frac{1}{\sqrt{n}}\le h_{1}$, for all $n\ge n_{0}$.

Then, for fixed arbitrary $n\ge n_{0}$, passing above to supremum after $h\in [0, 1/\sqrt{n})$, we immediately get
$$\omega_{\varphi}^{2}(f;1/\sqrt{n})\ge \lambda_{0}\cdot \frac{1}{n} \|\varphi^{2} f^{\prime \prime}\|,$$
which completes the proof. $\hfill\Box$

The crucial step in the proof of Theorem 1 is to establish norm estimate in the Theorem of Voronovskaja. This theorem was first proved in [20] and is given in the book of DeVore and Lorentz [1] as follows:

{\bf Theorem D.} {\it If $f$ is bounded on $[0,1]$, differentiable in some neighborhood of $x$ and has second derivative $f''(x)$ for some $x\in [0,1]$, then
$$\lim\limits_{n\to \infty} n\cdot[B_n(f,x)-f(x)]=\frac{x(1-x)}{2}\cdot f''(x).\eqno(2.6)$$
If $f\in C^2[0,1]$, the convergence is uniform.}

This result has attracted the attention of many authors in the last 80 years. Very recently some quantitative estimates in pointwise variant of Voronovskaja's theorem are obtained in [11,12,13,14,18] and for the complex-valued functions of complex variable in [4,5,6]. Concerning norm estimates we cite Lemma 8.3 in [2] as

{\bf Theorem E.} {\it For $f\in W^3_{\infty}(\varphi)$ and $n\ge 12$ we have }
$$
\|B_nf-f-\frac{1}{2n}\cdot \varphi^2f''\|_{L_{\infty}}\le n^{-\frac{3}{2}}\cdot \|\varphi^3f'''\|_{L_{\infty}}.\eqno(2.7)$$

Our next two statements extends this result to $f\in C^{2}[0, 1]$.

{\bf Theorem 3.} {\it For any $f\in C^{2}[0, 1]$, there exists $n_{0}:=n_{0}(f)$, such that for all $n\ge n_{0}$ we have }
$$
\left \|B_nf-f-\frac{1}{2n}\cdot \varphi^2f^{\prime \prime}\right \|\le 4 \omega^2_{\varphi}(f,\frac{1}{\sqrt{n}}).\eqno(2.8)$$

{\bf Proof:} Obviously
$$
\left \|B_nf-f-\frac{1}{2n}\cdot \varphi^2f^{\prime \prime}\right \|\le \|B_nf-f\|+\frac{1}{2n}\|\varphi^2f^{\prime \prime}\|\le $$
$$ \le 3\omega^2_{\varphi}(f,\frac{1}{\sqrt{n}})+\omega^2_{\varphi}(f,\frac{1}{\sqrt{n}})=4\omega^2_{\varphi}(f,\frac{1}{\sqrt{n}}),$$
where we have applied the upper estimate in [10] and (2.5), respectively. The proof is completed. $\hfill\Box$

Moreover, in terms of usual moduli of continuity, we can obtain the following better estimate.

{\bf Theorem 4.} {\it If $f\in C^{2}[0, 1]$, then for all $n\in \mathbb{N}$, $n\ge 2$, it holds
$$\|B_nf-f-\frac{1}{2n}\cdot \varphi^2f^{\prime \prime}\|\le \frac{5}{8 n}\omega_{1}(f^{\prime \prime}, \frac{1}{\sqrt{n}})+\frac{13}{64 n}\omega_{2}(f^{\prime \prime}, \frac{1}{\sqrt{n}}). \eqno(2.9)$$}

{\bf Proof.} By Theorem 4 in Gonska-Ra\c{s}a \cite{hg4}, for $f\in C^{2}[0, 1]$, $n\ge 2$, $x\in [0, 1]$ and $X=x(1-x)$, we have
$$\left |B_{n}(f, x)-f(x)-\frac{x(1-x)}{2n}f^{\prime \prime}(x)\right |$$
$$\le \frac{X}{n}\left \{\frac{X^{\prime}}{\sqrt{3(n-2)X+1}}\cdot \frac{5}{6}\omega_{1}\left (f^{\prime \prime}, \sqrt{\frac{3(n-2)X+1}{n^{2}}}\right )\right .$$
$$\left .+\frac{13}{16}\omega_{2}\left (f^{\prime \prime}, \sqrt{\frac{3(n-2)X+1}{n^{2}}}\right ) \right \}.$$
Passing to uniform norm and taking into account the inequalities $0\le X\le 1/4$, $|X^{\prime}|=|1-2x|\le 3$, for all $x\in [0, 1]$
we immediately obtain
$$\|B_nf-f-\frac{1}{2n}\cdot \varphi^2f^{\prime \prime}\|\le \frac{1}{4n}\left \{3\cdot \frac{5}{6}\omega_{1}\left (f^{\prime \prime},
\sqrt{\frac{[3(n-2)/4]+1}{n^{2}}}\right )\right .$$
$$\left .+\frac{13}{16}\omega_{2}\left (f^{\prime \prime},\sqrt{\frac{[3(n-2)/4]+1}{n^{2}}}\right )\right \}$$
$$\le \frac{5}{8n}\omega_{1}(f^{\prime \prime}, \frac{1}{\sqrt{n}})+\frac{13}{64 n}\omega_{2}(f^{\prime \prime}, \frac{1}{\sqrt{n}}),$$
which proves the theorem. $\hfill\Box$

\section{Proofs of Theorem 1 and Corollary 1}

{\bf Proof of Theorem 1.} Firstly,  note that in the statement of Theorem 1, we may suppose that $f$ is not a polynomial of degree $\le 1$, because if $f$ is a polynomial of degree $\le 1$ then Theorem 1 holds trivially. This supposition obviously implies that in what follows we have $\|\varphi^{2} f^{\prime \prime}\|>0$.

We apply the ideas in the case of Bernstein polynomials of complex variable in the proof of Theorem 2.1 in [8]. The following identity is valid for all $f$, which are not a polynomial of degree $\le 1$
$$
\begin{array}{l}
\ds B_n(f,x)-f(x)=\omega^2_{\varphi}(f,\frac{1}{\sqrt{n}})\Big\{\frac{\varphi^2(x)f''(x)}{2}\cdot \frac{1}{n\omega^2_{\varphi}(f,\frac{1}{\sqrt{n}})}+\\[6mm]
\ds + \left[\frac{1}{\omega^2_{\varphi}(f,\frac{1}{\sqrt{n}})}\left(B_n(f,x)-f(x)-\frac{\varphi^2(x)f''(x)}{2n}\right)\right]\Big\}.
\end{array}
\eqno(3.1)$$
This immediately implies
$$
\begin{array}{l}
\ds \|B_nf-f\|\ge \omega^2_{\varphi}(f,\frac{1}{\sqrt{n}})\Big\{\frac{\|\varphi^2f''\|}{2}\cdot \frac{1}{n\omega^2_{\varphi}(f,\frac{1}{\sqrt{n}})}-\\[6mm]
\ds - \left[\frac{1}{\omega^2_{\varphi}(f,\frac{1}{\sqrt{n}})}\cdot\|B_nf-f-\frac{\varphi^2f''}{2n}\|\right]\Big\}.
\end{array}
\eqno(3.2)$$
Now from (2.4) we obtain
$$\frac{1}{n\omega^2_{\varphi}(f,\frac{1}{\sqrt{n}})}\ge \frac{A_1}{\|\varphi^2f''\|}, \mbox{ for all } n\in \mathbb{N}, \eqno(3.3)$$
which implies
$$
\|B_nf-f\|\ge \omega^2_{\varphi}(f,\frac{1}{\sqrt{n}})\cdot \left\{\frac{A_1}{2}-A_{n}\right \}, \eqno(3.4)$$
where
$$A_n=\frac{\frac{5}{8n}\omega_{1}(f^{\prime \prime}, 1/\sqrt{n})+\frac{13}{64 n}\omega_{2}(f^{\prime \prime}, 1/\sqrt{n})}{\omega^2_{\varphi}(f,\frac{1}{\sqrt{n}})}\to 0,\,\mbox{as} \, n\to \infty.$$
More precisely, using (2.4) and (2.5), we immediately obtain the double inequality
$$\frac{1}{16\|\varphi^{2}f^{\prime \prime}\|}\left [\frac{5}{8}\omega_{1}(f^{\prime \prime} ; 1/\sqrt{n})+\frac{13}{64}\omega_{2}(f^{\prime \prime} ; 1/\sqrt{n})\right ]\le A_{n}\le$$
$$\le \frac{1}{\lambda_{0}\|\varphi^{2}f^{\prime \prime}\|}\left [\frac{5}{8}\omega_{1}(f^{\prime \prime} ; 1/\sqrt{n})+\frac{13}{64}\omega_{2}(f^{\prime \prime} ; 1/\sqrt{n})\right ], \eqno(3.5)$$
where the right-hand side holds for all $n\ge n_{0}(f, \lambda_{0})$ (that comes from $(2.5)$), while the left-hand side holds for all $n\ge 2$.

Therefore, from the right-hand side of (3.5) it follows that $\lim_{n\to \infty}A_{n}=0$.

From (3.4) and (3.5), evidently that there exists $n_{1}(f, \lambda_{0}, \mu_{0})>n_{0}(f, \lambda_{0})$, such that for all $n\ge n_{1}(f, \lambda_{0}, \mu_{0})$
$$\frac{A_{1}}{2}-A_{n}\ge \frac{\mu_{0} A_{1}}{2}=\frac{\mu_{0}}{32}.\eqno(3.6)$$
Therefore, from (3.4) we have
$$
\|B_nf-f\|\ge \frac{\mu_{0}}{32}\cdot \omega^2_{\varphi}(f,\frac{1}{\sqrt{n}}), \eqno(3.7)$$
for all $n\ge n_{1}$, which proves Theorem 1. $\hfill\Box$

{\bf Proof of Corollary 1.} We have to prove just the left-hand side inequality. Following the lines in the proof of Theorem 2, since
$\frac{1}{2}\in [h^{2}/(1+h^{2}), 1/(1+h^{2})]$ for all $0\le h <1$, we easily get
$$\sup  \{|f(x+h\varphi(x))-2f(x)+f(x-h\varphi(x))|; x\in [h^{2}/(1+h^{2}), 1/(1+h^{2})] \}$$
$$\ge |f(1/2+h\varphi(1/2)-2f(1/2)+f(1/2-h\varphi(1/2)|=h^{2}\varphi^{2}(1/2)\cdot |f^{\prime \prime}(\xi_{h, 1/2})|$$
$$\ge \frac{m h^{2}}{4}, \mbox{ for all } 0\le h<1,$$
where passing to supremum with $0\le h\le \delta$ implies $\omega_{\varphi}^{2}(f, \delta)\ge \frac{m}{4}\cdot \delta^{2}$.

Therefore, by taking $\delta=\frac{1}{\sqrt{n}}$, it follows
$\omega_{\varphi}^{2}(f, 1/\sqrt{n})\ge \frac{m}{4}\cdot \frac{1}{n}.$

Now, following the lines in the proof of Theorem 1, where instead of Theorem 4 we use Theorem E, we get the same relationship (3.4),
where because of the hypothesis, the estimate in Theorem E becomes
$$\|B_nf-f-\frac{1}{2n}\cdot \varphi^2f''\|\le \frac{M}{8n^{3/2}},$$
with $A_{n}$ upper bounded as follows :
$$A_{n}\le \frac{M n^{-3/2}}{8\omega_{\varphi}^{2}(f, 1/\sqrt{n})}\le \frac{M}{8n^{3/2}}\cdot \frac{4 n}{m}=\frac{M}{2m\sqrt{n}}.$$
Therefore, from the lines in the proof of Theorem 1 we obtain
$$\|B_nf-f\|\ge \omega^2_{\varphi}(f,\frac{1}{\sqrt{n}})\cdot \left\{\frac{A_1}{2}-\frac{M}{2m\sqrt{n}}\right \},$$
which by (3.6) is valid for all $n$ satisfying
$$\frac{M}{2m\sqrt{n}}\le \frac{1}{64},$$
that is for all $n\ge n_{1}$ with $n_{1}=\left [\frac{1024 M^{2}}{m^{2}}\right ]+1$, which proves the corollary. $\hfill\Box$

{\bf Remark 3.} If we suppose, in addition, that for example $f\in C^{4}[0, 1]$, then from the above proof of Theorem 1, relationship (3.6), it easily follows that $A_{n}\le \frac{1-\mu_{0}}{32}$, and therefore that the index $n_{1}$ in the statement, necessarily must be chosen greater than the smallest number $n_{2}\in \mathbb{N}$ that satisfies the inequality
$$\frac{1}{\lambda_{0}\|\varphi^{2}f^{\prime \prime}\|}\left [\frac{5}{8\sqrt{n}}\cdot \|f^{\prime \prime \prime}\|+\frac{13}{64 n}\cdot
\|f^{(4)}\|\right ]\le \frac{1-\mu_{0}}{32}.$$
Indeed, this immediately follows from the inequality
$$\frac{1}{\lambda_{0}\|\varphi^{2}f^{\prime \prime}\|}\left [\frac{5}{8}\omega_{1}(f^{\prime \prime} ; 1/\sqrt{n})+\frac{13}{64}\omega_{2}(f^{\prime \prime} ; 1/\sqrt{n})\right ]$$
$$\le \frac{1}{\lambda_{0}\|\varphi^{2}f^{\prime \prime}\|}\left [\frac{5}{8\sqrt{n}}\cdot \|f^{\prime \prime \prime}\|+\frac{13}{64 n}\cdot\|f^{(4)}\|\right ].$$

If instead of $f\in C^{4}[0, 1]$, we suppose that $f\in W^{3}_{\infty}(\varphi)$ and in the proof of Theorem 1, instead of Theorem 4 we use Theorem E, then we easily get
$$A_{n}=\frac{n^{-3/2}\|\varphi^{3}f^{\prime \prime \prime}\|}{[\lambda_{0}\|\varphi^{2} f^{\prime \prime}\|/(n)]}
=\frac{1}{\sqrt{n}}\cdot \frac{\|\varphi^{3} f^{\prime \prime \prime}\|}{\lambda_{0} \|\varphi^{2} f^{\prime \prime}\|},$$
for all $n\ge n_{0}(f, \lambda_{0})$ and that the index $n_{1}(f,\mu_{0}, \lambda_{0})$ in the statement of Theorem 1, necessarily must be greater than the smallest $n_{2}$ that satisfies the inequality $A_{n}\le \frac{1-\mu_{0}}{32}$. Simple calculation shows that we may take
$$n_{2}=\left [\left (\frac{32 \|\varphi^{3} f^{\prime \prime \prime}\|}{\lambda_{0}(1-\mu_{0})\|\varphi^{2} f^{\prime \prime}\|}\right )^{2}\right ]+1,$$
where $[a]$ means the integer part of $a$.

{\bf Remark 4.} Note that for the limit case of Bernstein operator, $U_n(f,x)$, it was proved by Parvanov and Popov the following strong converse inequality in [16]:
$$\frac{1}{2}\|U_nf-f\|\le K^2_{\varphi}(f,\frac{1}{n})\le (6+\sqrt{8})\|U_nf-f\|.\eqno(3.8)$$
The proof relies on the commutativity of $U_n$-a property, which is not available for $B_n$.

\section{Concrete Examples}

For some particular classes of functions, the index $n_{1}$ in Theorem 1 and Corollary 1 can be explicitly obtained,
as follows.

{\bf Example 1.} Firstly, consider $f(x)=exp(x)$, $x\in [0, 1]$. Since we have $\|f^{\prime \prime \prime}\|\le e$ and $|f^{\prime \prime}(x)|\ge 1$, for all $x\in[0,1]$, it follows that we can apply Corollary 1, obtaining that we have
$$\|B_{n}(exp, \cdot)-exp\|\ge \frac{1}{64}\omega_{\varphi}^{2}(exp, 1/\sqrt{n}), \mbox{ for all } n> 1024 e^{2}.$$

{\bf Example 2.} Secondly,consider $f(x)=cos(x)$, $x\in [0, 1]$. Since have  $\|f^{\prime \prime \prime}\|\le 1$ and $|f^{\prime \prime}(x)|\ge cos(1)\approx 0.540302306$, for all $x\in[0,1]$, it follows
that again we can apply Corollary 1, obtaining that we have
$$\|B_{n}(cos, \cdot)-cos\|\ge \frac{1}{64}\omega_{\varphi}^{2}(cos, 1/\sqrt{n}), \mbox{ for all } n\ge \left [\frac{1024}{cos(1)}\right ]+1.$$

{\bf Example 3.} Thirdly, by taking $f(x)=sin(x)$, $x\in [0, 1]$, it is clear that we cannot apply Corollary 1, but we can apply Theorem 1, or more exactly its variant expressed by Remark 3, first part.

Firstly, let us find the index $n_{0}$ in the Theorem 2 (as Theorem 2 is used in the proof of Theorem 1). For this purpose,
note that simple calculation leads us to
$$|f(x+h\varphi(x))-2f(x)+f(x-h\varphi(x))|=2|sin(x)[1-cos(h\varphi(x))]|$$
$$=4sin(x)\cdot sin^{2}(h\varphi(x)/2)\ge 4\cdot \frac{2}{\pi}x\cdot \frac{4}{\pi^{2}}\cdot \frac{h^{2}\varphi^{2}(x)}{4}=\frac{8}{\pi^{3}}x h^{2}\varphi^{2}(x),$$
for all $x\in \left [\frac{h^{2}}{1+h^{2}}, \frac{1}{1+h^{2}}\right ]$. Here we used the inequality $sin(x)\ge \frac{2}{\pi} x$, for all
$x\in [0, \pi/2]$.

Suppose in what follows that $0\le h\le \delta \le \frac{1}{\sqrt{2}}$. It follows that $\frac{2}{3}\in \left [\frac{h^{2}}{1+h^{2}}, \frac{1}{1+h^{2}}\right ]$. Also, clearly that $x\varphi^{2}(x)$ attains its maximum value $\frac{4}{27}$ on $[0, 1]$ at $x=2/3$.

Therefore, passing above to supremum firstly with respect to $x$, and then with respect to $h\in [0, \delta]$, with $\delta\le \frac{1}{\sqrt{2}}$,
it is immediate that
$$\omega_{\varphi}^{2}(sin, \delta)\ge \frac{8}{\pi^{3}}\cdot \frac{4}{27}\cdot \delta^{2}=\frac{32}{27 \pi^{3}}\delta^{2},$$
that is for $\delta=1/\sqrt{n}$ with $n\ge 2$, we obtain
$$\omega_{\varphi}^{2}(sin, 1/\sqrt{n})\ge \frac{32}{27 \pi^{3}}\cdot \frac{1}{n}\ge \frac{32}{27 \pi^{3}}\cdot \frac{\|\varphi^{2} f^{\prime \prime}\|}{n}, \mbox{ for all } n\ge 2.$$

It follows that Theorem 2 holds with $\lambda_{0}=\frac{32}{27 \pi^{3}}$ and $n_{0}=2$.

Now, because
$$\|\varphi^{2} f^{\prime \prime}\|\ge \varphi^{2}(1/2) sin(1/2)=\frac{1}{4}sin(1/2)\approx 0.25 \times 0.4794=0.11985,$$
by the first part of Remark 3 we obtain that $n_{2}$ there, can be chosen as the smallest index $n$ satisfying the last inequality below
$$\frac{1}{\lambda_{0}\|\varphi^{2}f^{\prime \prime}\|}\left [\frac{5}{8\sqrt{n}}\cdot \|f^{\prime \prime \prime}\|+\frac{13}{64 n}\cdot
\|f^{(4)}\|\right ]\le \frac{4}{\lambda_{0}sin(1/2)}\left [\frac{5}{8 \sqrt{n}}+\frac{13}{64 n}\right ]$$
$$\le \frac{1-\mu_{0}}{32}, \eqno{(4.1)}$$
with $\lambda_{0}, \mu_{0}\in (0, 1)$, arbitrary fixed.

From here, it is immediate that $n_{2}$ can be chosen the smallest index $n$ satisfying
$$\frac{4}{\lambda_{0}sin(1/2)}\left [\frac{5}{8}+\frac{13}{64}\right ]\cdot \frac{1}{\sqrt{n}}\le \frac{1-\mu_{0}}{32},$$
which by simple calculation leads to $n_{2}=\left [\left (\frac{212}{\lambda_{0}(1-\mu_{0})sin(1/2)}\right )^{2}\right ]+1$
(here $[a]$ means the integer part of $a$), with $\lambda_{0}=\frac{32}{27 \pi^{3}}$. Choosing above $\mu_{0}=1/2$, we obtain
$$\|B_{n}(sin, \cdot)-sin\|\ge \frac{1}{64}\omega^{2}_{\varphi}(sin, \frac{1}{\sqrt{n}}), \mbox{ for all } n\ge \max\{2, n_{2}\}=n_{2}.$$

$\begin{array}{ll}
\textrm{S. G. Gal} &\textrm{G. T. Tachev} \\
 \textrm{University of Oradea} &\textrm{University of Architecture} \\
\textrm{Dept. of Math. and Comp. Sci.} &\textrm{ Dept. of Math.}\\
\textrm{RO-410087 Oradea} &\textrm{ BG-1046 Sofia}\\
\textrm{Romania}&\textrm{ Bulgaria}\\
\textrm{e-mail: galso@uoradea.ro}& \textrm{e-mail: gtt\_fte@uacg.bg}\\
\end{array}
$

\end{document}